\documentclass[10pt]{amsart}
\usepackage{amsfonts}
\usepackage{ifthen}
\usepackage{amsthm}
\usepackage{amsmath}
\usepackage{graphicx}
\usepackage{amscd,amssymb,amsthm}
\usepackage{graphicx}

\newcounter{minutes}
\divide\time by 60
\newcounter{hours}
\multiply\time by 60 \addtocounter{minutes}{-\time}

\title[The radius of starlikeness of normalized Bessel
functions]{The radius of starlikeness of normalized Bessel
functions of the first kind}

\author[\'A. Baricz]{\'Arp\'ad Baricz}
\address{Department of Economics, Babe\c{s}-Bolyai University,
Cluj-Napoca 400591, Romania} \email{arpad.baricz@econ.ubbcluj.ro}

\author[P.A. Kup\'an]{P\'al Aurel Kup\'an}
\address{Department of Mathematics, Sapientia University, T\^irgu
Mure\c{s} 540485, Romania} \email{kupanp@ms.sapientia.ro}

\author[R. Sz\'asz]{R\'obert Sz\'asz}
\address{Department of Mathematics, Sapientia University, T\^irgu
Mure\c{s} 540485, Romania} \email{rszasz@ms.sapientia.ro}

\thanks{\'A. Baricz was supported by the J\'anos
Bolyai Research Scholarship of the Hungarian Academy of Sciences. P.A. Kup\'an and R. Sz\'asz were supported by the
Sapientia Research Foundation.}

\newtheorem{theorem}{Theorem}
\newtheorem{corollary}{Corollary}

\newcommand{\real}{\operatorname{Re}}

\newcommand{\loga}{\operatorname{Log}}
\keywords{Bessel and modified Bessel functions of the first kind;
univalent, starlike and convex functions.} \subjclass[2010]{30C45,
33C10}
\begin{document}

\def\thefootnote{}
\footnotetext{ \texttt{File:~\jobname .tex,
          printed: \number\year-0\number\month-0\number\day,
          \thehours.\ifnum\theminutes<10{0}\fi\theminutes}
} \makeatletter\def\thefootnote{\@arabic\c@footnote}\makeatother

\maketitle
%=======================================================================================================================================================

\begin{abstract}
In this note our aim is to determine the radius of starlikeness of
the normalized Bessel functions of the first kind for three
different kinds of normalization. The key tool in the proof of our
main result is the Mittag-Leffler expansion for Bessel functions of
the first kind and the fact that, according to Ismail and Muldoon \cite{ismail}, the smallest positive zeros of
some Dini functions are less than the first positive zero of the
Bessel function of the first kind.
\end{abstract}

%=======================================================================================================================================================

\section{\bf Introduction}

Let $\mathbf{D}(z_0,r) = \{z\in \mathbb{C}: |z-z_0| < r\}$ be the
open disk with center $z_{0}\in\mathbb{C}$ and radius $r>0$ and let
us denote the particular disk $\mathbf{D}(0,1)$ by $\mathbf{D}.$
Moreover, let $\mathcal{A}$ be the class of analytic functions
defined in the unit disk $\mathbf{D}$, which can be normalized as
$f(z)=z+a_{2}z^{2}+\dots,$ that is, $f(0)=f'(0)-1=0.$ The class of
starlike functions, denoted by $\mathcal{S}^*,$ is the subclass of
$\mathcal{A}$ which consists of functions $f$ for which the domain
$f(\mathbf{D})$ is starlike with respect to $0.$ An analytic
description of $\mathcal{S}^{*}$ is
$$\mathcal{S}^* = \left\{f\in \mathcal{A}\left|
\real\left[\frac{zf'(z)}{f(z)}\right]\right. > 0 \ \ \mbox{for all}
\ \ z\in \mathbf{D} \right\}.$$ Moreover, consider the class of starlike
functions of order $\beta\in[0,1),$ that is,
$$\mathcal{S}^*(\beta) = \left\{f\in \mathcal{A}\left|
\real\left[\frac{zf'(z)}{f(z)}\right]\right. > \beta \ \ \mbox{for all}
\ \ z\in \mathbf{D} \right\}.$$
The real numbers
$$r^{*}(f)=\sup\left\{r>0\left| \real\left[\frac{zf'(z)}{f(z)}\right]\right. > 0\ \ \mbox{for all} \ \
z\in \mathbf{D}(0,r)\right\}$$
and
$$r^{*}_{\beta}(f)=\sup\left\{r>0\left| \real\left[\frac{zf'(z)}{f(z)}\right]\right. > \beta\ \ \mbox{for all} \ \
z\in \mathbf{D}(0,r)\right\},$$
are called the radius of starlikeness and the radius of
starlikeness of order $\beta$ of the function $f,$ respectively. We note that in fact $r^{*}(f)$ is the largest radius such that
$f\left(\mathbf{D}(0,r^{*}(f))\right)$ is a starlike domain with
respect to $0.$

Now, consider the Bessel function of the first kind \cite{watson}, which is a particular solution of the
second-order linear homogeneous Bessel differential equation. This function has the
infinite series representation
$$J_{\nu}(z)=\sum_{n\geq0}\frac{(-1)^{n}}{n!\Gamma(n+\nu+1)}\left(\frac{z}{2}\right)^{2n+\nu},$$
where $z\in\mathbb{C}$ and $\nu\in\mathbb{C}$ such that $\nu\neq-1,-2,{\dots}.$
Observe that the Bessel function $J_{\nu}$ does not belong to class
$\mathcal{A}.$ Thus, it is natural to consider the following three
kind of normalization of the Bessel function of the first kind
\begin{equation}\label{f}f_{\nu}(z)=\left[2^{\nu}\Gamma(\nu+1)J_{\nu}(z)\right]^{1/\nu},
\nu\neq{0},\end{equation}
\begin{equation}\label{g}g_{\nu}(z)=2^{\nu}\Gamma(\nu+1)z^{1-\nu}J_{\nu}(z)\end{equation}
and
\begin{equation}\label{h}h_{\nu}(z)=2^{\nu}\Gamma(\nu+1)z^{1-{\nu}/{2}}J_{\nu}(\sqrt{z}).\end{equation}
Clearly the functions $f_{\nu},$ $g_{\nu}$ and $h_{\nu}$ belong to the class
$\mathcal{A}.$ We note here that in fact
$$f_{\nu}(z)=\exp\left[\frac{1}{\nu}\loga\left(2^{\nu}\Gamma(\nu+1)J_{\nu}(z)\right)\right],$$
where $\loga$ represents the principal branch of the logarithm
function and every many-valued function considered in this paper are
taken with the principal branch.

Now, let us recall some results on the geometric behavior of the functions $f_{\nu},$ $g_{\nu}$ and $h_{\nu}.$ Brown \cite{brown} determined
the radius of starlikeness for $f_{\nu}$ in the case when $\nu>0.$ Namely, in \cite[Theorem 2]{brown} it was shown that the radius $r^{*}\left(f_{\nu}\right)$
is the smallest positive zero of the function $z\mapsto J_{\nu}'(z).$ Moreover, in \cite[Theorem 3]{brown} Brown proved that if $\nu>0,$ then
the radius of starlikeness of the function $g_{\nu}$ is the smallest positive zero of the function $z\mapsto zJ_{\nu}'(z)+(1-\nu)J_{\nu}(z).$
Kreyszig and Todd \cite[Theorem 3]{todd} proved that when $\nu>-1$ the function $g_{\nu}$ is univalent in the circle $|z|\leq\rho_\nu$ but not in any concentric circle with larger radius, where $\rho_\nu$ is the first maximum of the function $g_{\nu}$ on the positive real axis. Brown \cite[p. 282]{brown} pointed out that when $\nu>0$ the radius of starlikeness of the function $g_{\nu},$ that is, $r^{*}\left(g_{\nu}\right)$  is exactly the radius of univalence $\rho_{\nu}$ obtained by Kreyszig and Todd \cite{todd}.
Furthermore, Brown \cite[Theorem 5.1]{brown2} showed that the radius of starlikeness of the function $g_{\nu}$ is also $\rho_{\nu}$ when $\nu\in(-1/2,0).$ On the other hand, Hayden and Merkes \cite[Theorem C]{hayden} deduced that when $\mu=\real \nu>-1$ the radius of starlikeness of $g_{\nu}$ is not less than the smallest positive zero of $g_{\mu}'.$ It is worth to mention that Brown used the methods of Nehari \cite{nehari} and Robertson \cite{robertson}, and an important tool in the proofs was the fact that the Bessel function of the first kind is a particular solution of the Bessel differential equation. For related (more general) results the interested reader is referred to \cite{brown3,merkes,robertson,wilf} and to the references therein. Finally, let us mention that other geometric properties of the functions $g_{\nu}$ and $h_{\nu}$ were obtained in \cite{baricz1,bariczbook,bariczpon,szasz1,szasz2}. See also the references therein.

Motivated by the above results in this paper we make a contribution to the subject and we determine the radius of starlikeness of order $\beta$ for the functions $f_{\nu},$ $g_{\nu}$ and $h_{\nu}.$ We note that our approach is much simpler than the methods used in \cite{brown,brown2,hayden,todd}, and is based only on the Mittag-Leffler expansion for Bessel functions of the first kind and on the fact that the smallest positive zeros of certain Dini functions are less than the first positive zero of the Bessel function of the first kind, according to Ismail and Muldoon \cite{muldoon,ismail}.

\section{\bf Starlikeness of order $\beta$ of normalized Bessel functions}
\setcounter{equation}{0}

Our main result is the following theorem. Here $I_{\nu}$ denotes the modified Bessel function of the first kind, which in view of the relation $I_{\nu}(z)=\mathrm{i}^{-\nu}J_{\nu}(\mathrm{i}z)$ is also called sometimes as the Bessel function of the first kind with imaginary argument.

\begin{theorem}\label{th}
Let $1>\beta\geq0.$ Then the following assertions are true:
\begin{enumerate}
\item[\bf a.] If $\nu\in(-1,0),$ then $r^{*}_{\beta}\left(f_{\nu}\right)=x_{\nu,\beta},$ where $x_{\nu,\beta}$ denotes the unique positive root of the equation $zI_{\nu}'(z)-\beta\nu{I_{\nu}}(z)=0.$ Moreover, if $\nu>0,$ then we have $r^{*}_{\beta}\left(f_{\nu}\right)=x_{\nu,\beta,1},$ where $x_{\nu,\beta,1}$ is the smallest positive root of the equation $zJ_{\nu}'(z)-\beta\nu{J_{\nu}}(z)=0.$
\item[\bf b.] If $\nu>-1,$ then $r^{*}_{\beta}\left(g_{\nu}\right)=y_{\nu,\beta,1},$ where $y_{\nu,\beta,1}$ is the smallest positive root of the equation $zJ_{\nu}'(z)+(1-\beta-\nu){J_{\nu}}(z)=0.$
\item[\bf c.] If $\nu>-1,$ then $r^{*}_{\beta}\left(h_{\nu}\right)=z_{\nu,\beta,1},$ where $z_{\nu,\beta,1}$ is the smallest positive root of the equation $zJ_{\nu}'(z)+(2-2\beta-\nu){J_{\nu}}(z)=0.$
\end{enumerate}
\end{theorem}

In particular, when $\beta=0,$ we get the following result.

\begin{corollary}\label{cr}
The following assertions are true:
\begin{enumerate}
\item[\bf a.] If $\nu\in(-1,0),$ then the radius of starlikeness of $f_{\nu}$ is $x_{\nu,0},$ where $x_{\nu,0}$ is the unique positive root of the equation $I_{\nu}'(z)=0.$ If $\nu>0,$ then the radius of starlikeness of the function $f_{\nu}$ is $x_{\nu,0,1},$ which denotes the smallest positive root of the equation $J_{\nu}'(z)=0.$
\item[\bf b.] If $\nu>-1,$ then the radius of starlikeness of the function $g_{\nu}$ is $y_{\nu,0,1},$ which denotes the smallest positive root of the equation $zJ_{\nu}'(z)+(1-\nu){J_{\nu}}(z)=0.$
\item[\bf c.] If $\nu>-1,$ then the radius of starlikeness of the function $h_{\nu}$ is $z_{\nu,0,1},$ which denotes the smallest positive root of the equation $zJ_{\nu}'(z)+(2-\nu){J_{\nu}}(z)=0.$
\end{enumerate}
\end{corollary}

Observe that part {\bf a} and {\bf b} of Corollary \ref{cr} complement the results of \cite[Theorem 2]{brown}, \cite[Theorem 3]{brown} and \cite[Theorem 5.1]{brown2}, mentioned above. Part {\bf c} complements the results from \cite{baricz1,bariczpon,szasz1,szasz2}. It is of interest to note here that very recently Sz\'asz \cite{szasz1} proved that the normalized Bessel function $h_{\nu}$ is starlike if and only if $\nu\geq\nu_0,$ where $\nu_0=-0.5623\dots$ is the root of the equation $h_{\nu}'(1)=0,$ that is, $J_{\nu}'(1)+(2-\nu)J_{\nu}(1)=0.$ Finally, we mention that if we consider the function $z\mapsto\lambda_{\nu}(z)=h_{\nu}(z)/z,$ then part {\bf c} of Theorem \ref{th} and Corollary \ref{cr} can be rewritten in terms of convex functions. The idea is to use the differentiation formula
$$\lambda_{\nu}'(z)=-\frac{1}{4(\nu+1)}\lambda_{\nu+1}(z)$$
together with the well-known duality theorems of Alexander \cite{alexander} and Jack \cite{jack}. See also \cite[p. 25]{bariczbook} for the results of Alexander and Jack, and also \cite[Ch. 2]{bariczbook} for similar results on convex Bessel functions.

\begin{proof}[\bf Proof of Theorem \ref{th}]
First we prove part {\bf a} for $\nu>0$ and parts {\bf b} and {\bf c} for $\nu>-1.$ We need to show that the inequalities
\begin{equation}\label{starbeta}\real\left[\frac{zf_{\nu}'(z)}{f_{\nu}(z)}\right]>\beta,\ \ \real\left[\frac{zg_{\nu}'(z)}{g_{\nu}(z)}\right]>\beta\ \ \ \mbox{and}\ \ \ \real\left[\frac{zh_{\nu}'(z)}{h_{\nu}(z)}\right]>\beta\end{equation}
are valid for all $\nu>0$ and $z\in\mathbf{D}(0,x_{\nu,\beta,1}),$
$\nu>-1$ and $z\in\mathbf{D}(0,y_{\nu,\beta,1}),$ and $\nu>-1$ and $z\in\mathbf{D}(0,z_{\nu,\beta,1}),$ respectively, and each of the above inequalities does not hold in any larger disk.

Lommel's well-known result states that if $\nu>-1,$ then the
zeros of the Bessel function $J_{\nu}$ are all real. Thus, if
$j_{\nu,n}$ denotes the $n$-th positive zero of the Bessel function
$J_{\nu},$ then the Bessel function admits the Weierstrassian
decomposition of the form \cite[p. 498]{watson}
\begin{equation}\label{prod}J_{\nu}(z)=\frac{z^{\nu}}{2^{\nu}\Gamma(\nu+1)}\prod_{n\geq1}\left(1-\frac{z^2}{j_{\nu,n}^2}\right),\end{equation}
and this infinite product is uniformly convergent on each compact
subset of $\mathbb{C}.$ Logarithmic differentiation of \eqref{prod}
yields
\begin{equation}\label{logf}\frac{zJ'_{\nu}(z)}{J_{\nu}(z)}=\nu-\sum_{n\geq1}\frac{2z^{2}}{j_{\nu,n}^{2}-z^{2}},\end{equation}
which in view of the recurrence relation \cite[p. 45]{watson} $zJ_{\nu}'(z)-\nu
J_{\nu}(z)=-zJ_{\nu+1}(z)$ is equivalent to the Mittag-Leffler
expansion \cite[p. 498]{watson}
$$\frac{J_{\nu+1}(z)}{J_{\nu}(z)}=\sum_{n\geq1}\frac{2z}{j_{\nu,n}^{2}-z^{2}}.$$
Consequently, in view of \eqref{f},\eqref{g}, \eqref{h} and \eqref{logf} we obtain
$$\frac{zf'_{\nu}(z)}{f_{\nu}(z)}=\frac{1}{\nu}\frac{zJ'_{\nu}(z)}{J_{\nu}(z)}=1-\frac{1}{\nu}\sum_{n\geq1}\frac{2z^{2}}{j_{\nu,n}^{2}-z^{2}},$$
$$\frac{zg'_{\nu}(z)}{g_{\nu}(z)}=1-\nu+\frac{zJ'_{\nu}(z)}{J_{\nu}(z)}=1-\sum_{n\geq1}\frac{2z^{2}}{j_{\nu,n}^{2}-z^{2}}$$
and
$$\frac{zh'_{\nu}(z)}{h_{\nu}(z)}=1-\frac{\nu}{2}+\frac{1}{2}\frac{\sqrt{z}J'_{\nu}(\sqrt{z})}{J_{\nu}(\sqrt{z})}=1-\sum_{n\geq1}\frac{z}{j_{\nu,n}^{2}-z}.$$
It is known \cite[p. 597]{watson} that in case $\alpha+\nu>0$
and $\nu>-1$ the so-called Dini function $z\mapsto
zJ_{\nu}'(z)+\alpha J_{\nu}(z)$ has only real zeros and according to
Ismail and Muldoon \cite[p. 11]{ismail} we know that the smallest
positive zero of the above function is less than $j_{\nu,1}.$ This in turn implies
that $x_{\nu,\beta,1}<j_{\nu,1}$ for all $\nu>0,$ $y_{\nu,\beta,1}<j_{\nu,1}$ for all $\nu>-1,$ and
$z_{\nu,\beta,1}<j_{\nu,1}$ for all $\nu>-1.$ In other words, for all $\beta<1$ and $n\in\{2,3,{\dots}\}$ we have
$\mathbf{D}(0,x_{\nu,\beta,1})\subset\mathbf{D}(0,j_{\nu,1})\subset\mathbf{D}(0,j_{\nu,n})$
when $\nu>0,$ $\mathbf{D}(0,y_{\nu,\beta,1})\subset\mathbf{D}(0,j_{\nu,1})\subset\mathbf{D}(0,j_{\nu,n})$
when $\nu>-1,$ and $\mathbf{D}(0,y_{\nu,\beta,1})\subset\mathbf{D}(0,j_{\nu,1})\subset\mathbf{D}(0,j_{\nu,n})$
when $\nu>-1.$ On the other hand, it is known \cite{szasz1} that if $z\in{\mathbb{C}}$ and $\alpha\in{\mathbb{R}}$ such that
$\alpha>|z|$, then
\begin{equation}\label{a}\frac{|z|}{\alpha-|z|}\geq \real\left(\frac{z}{\alpha-z}\right).\end{equation}
By using \eqref{a}, we
obtain for all $\nu>-1,$ $n\in\{1,2,\dots\}$ and
$z\in\mathbf{D}(0,j_{\nu,1})$ the inequality
\begin{equation}\label{realj}\frac{|z|^{2}}{j_{\nu,n}^{2}-|z|^{2}}\geq
\real\left(\frac{z^{2}}{j_{\nu,n}^{2}-z^{2}}\right),\end{equation}
which in turn implies that
$$\real\left[\frac{zf'_{\nu}(z)}{f_{\nu}(z)}\right]=1-\frac{1}{\nu}\real\left[\sum_{n\geq1}\frac{2z^{2}}{j_{\nu,n}^{2}-z^{2}}\right]\geq
1-\frac{1}{\nu}\sum_{n\geq1}\frac{2|z|^{2}}{j_{\nu,n}^{2}-|z|^{2}}=\frac{|z|f'_{\nu}(|z|)}{f_{\nu}(|z|)},$$
$$\real\left[\frac{zg'_{\nu}(z)}{g_{\nu}(z)}\right]=1-\real\left[\sum_{n\geq1}\frac{2z^{2}}{j_{\nu,n}^{2}-z^{2}}\right]\geq
1-\sum_{n\geq1}\frac{2|z|^{2}}{j_{\nu,n}^{2}-|z|^{2}}=\frac{|z|g'_{\nu}(|z|)}{g_{\nu}(|z|)}
$$
and
$$\real\left[\frac{zh'_{\nu}(z)}{h_{\nu}(z)}\right]=1-\real\left[\sum_{n\geq1}\frac{z}{j_{\nu,n}^{2}-z}\right]\geq
1-\sum_{n\geq1}\frac{|z|}{j_{\nu,n}^{2}-|z|}=\frac{|z|h'_{\nu}(|z|)}{h_{\nu}(|z|)},$$
with equality when $z=|z|=r$. The minimum principle for harmonic
functions and the previous inequalities imply that the corresponding inequalities in \eqref{starbeta} are
valid if and only if we have $|z|<x_{\nu,\beta,1},$ $|z|<y_{\nu,\beta,1},$ and $|z|<z_{\nu,\beta,1},$ respectively, where $x_{\nu,\beta,1},$ $y_{\nu,\beta,1}$ and
$z_{\nu,\beta,1}$ are the smallest positive roots of the equations
$${rf'_{\nu}(r)}/{f_{\nu}(r)}=\beta,\ \ {rg'_{\nu}(r)}/{g_{\nu}(r)}=\beta$$
and
$${rh'_{\nu}(r)}/{h_{\nu}(r)}=\beta,$$ which are equivalent to
$$rJ_{\nu}'(r)-\beta\nu{J_{\nu}}(r)=0,\ \ rJ_{\nu}'(r)+(1-\beta-\nu){J_{\nu}}(r)=0$$ and $$rJ_{\nu}'(r)+(2-2\beta-\nu){J_{\nu}}(r)=0,$$ respectively.

Now, we prove the statement of part {\bf a} when $\nu\in(-1,0).$ First observe that the counterpart of \eqref{a}, that is,
\begin{equation}\label{aa}\real\left(\frac{z}{\alpha-z}\right)\geq \frac{-|z|}{\alpha+|z|}\end{equation}
is valid for all $\alpha\in\mathbb{R}$ and $z\in\mathbb{C}$ such that $\alpha>|z|.$ Indeed, if we have $z=x+\mathrm{i}y$ and $m=|z|=\sqrt{x^2+y^2},$
then \eqref{aa} is equivalent to $\alpha(\alpha-m)(m+x)\geq0,$ which is clearly true. By using \eqref{aa}, we
obtain for all $\nu>-1,$ $n\in\{1,2,\dots\}$ and
$z\in\mathbf{D}(0,j_{\nu,1})$ the inequality
\begin{equation}\label{realjj}
\real\left(\frac{z^{2}}{j_{\nu,n}^{2}-z^{2}}\right)\geq\frac{-|z|^{2}}{j_{\nu,n}^{2}+|z|^{2}},\end{equation}
which in turn implies that
$$\real\left[\frac{zf'_{\nu}(z)}{f_{\nu}(z)}\right]=1-\frac{1}{\nu}\real\left[\sum_{n\geq1}\frac{2z^{2}}{j_{\nu,n}^{2}-z^{2}}\right]\geq
1+\frac{1}{\nu}\sum_{n\geq1}\frac{2|z|^{2}}{j_{\nu,n}^{2}+|z|^{2}}=\frac{\mathrm{i}|z|f'_{\nu}(\mathrm{i}|z|)}{f_{\nu}(\mathrm{i}|z|)}.$$
This time we have equality if $z=\mathrm{i}|z|=\mathrm{i}r,$ and from the above inequality we conclude
that the first inequality in \eqref{starbeta} holds if and only if $|z|<x_{\nu,\beta},$ where
$x_{\nu,\beta}$ denotes the positive root of the equation ${\mathrm{i}rf'_{\nu}(\mathrm{i}r)}/{f_{\nu}(\mathrm{i}r)}=\beta,$ which is
equivalent to $rI_{\nu}'(r)-\beta\nu{I_{\nu}}(r)=0.$ All we need to prove is that $x_{\nu,\beta}$ is unique and $x_{\nu,\beta}<j_{\nu,1}$ for all $\beta\in[0,1)$ and $\nu\in(-1,0),$ since in order to use \eqref{realjj} we tacitly assumed that for all $\beta\in[0,1)$ and $n\in\{2,3,{\dots}\}$ we have
$\mathbf{D}(0,x_{\nu,\beta})\subset\mathbf{D}(0,j_{\nu,1})\subset\mathbf{D}(0,j_{\nu,n})$
when $\nu\in(-1,0).$ For this recall that in case $-1<\nu<-\alpha$ the Dini function $z\mapsto
zJ_{\nu}'(z)+\alpha J_{\nu}(z)$ has all its zeros real and a single pair of conjugate purely
imaginary zeros \cite[p. 597]{watson}. Moreover, due to
Ismail and Muldoon \cite[eq. (3.2)]{muldoon} we know that if $\pm\mathrm{i}\xi$ ($\xi$ real) denote the purely imaginary zeros of the Dini function $z\mapsto zJ_{\nu}'(z)+\alpha J_{\nu}(z),$ then
$$\xi^2<-\frac{\alpha+\nu}{2+\alpha+\nu}j_{\nu,1}^2.$$
This in turn implies that
$$x_{\nu,\beta}^2<-\frac{\nu(1-\beta)}{2+\nu(1-\beta)}j_{\nu,1}^2<j_{\nu,1}^2,$$
as we required. Finally, consider the function $q_{\nu}:(0,\infty)\to\mathbb{R},$ defined by $q_{\nu}(r)=rI_{\nu}'(r)/I_{\nu}(r)-\beta\nu.$ By using the asymptotic relations for small and large values of $r$ for the function $r\mapsto I_{\nu}(r),$ it can be verified that $rI_{\nu}'(r)/I_{\nu}(r)$ tends to $\nu$ as $r\to0,$ and tends to infinity as $r\to\infty.$ Moreover, it is known (see for example \cite{bariczbull}) that the function $r\mapsto rI_{\nu}'(r)/I_{\nu}(r)$ is increasing on $(0,\infty)$ for all $\nu>-1.$ Thus the function $q_{\nu}$ is increasing, $q_{\nu}(r)$ tends to $\nu(1-\beta)<0$ as $r\to0,$ and tends to infinity as $r\to\infty.$ Consequently, the graph of $q_{\nu}$ intersects the $r$-axis only once, and thus the equation $rI_{\nu}'(r)-\beta\nu{I_{\nu}}(r)=0$ has only one solution. This completes the proof.
\end{proof}


\begin{thebibliography}{width}

\bibitem[\bf Al]{alexander} \textsc{J.W. Alexander}, Functions which map the interior of the unit circle upon simple regions,
{\em Ann. of Math.} 17 (1915) 12--29.

\bibitem[\bf Ba1]{baricz1} \textsc{\'A. Baricz}, Geometric properties of generalized Bessel functions, {\em Publ.
Math. Debrecen} 73 (2008) 155--178.

\bibitem[\bf Ba2]{bariczbull} \textsc{\'A. Baricz}, Tur\'an type inequalities for modified Bessel functions, {\em Bull. Aust. Math. Soc.} 82 (2010) 254--264.

\bibitem[\bf Ba3]{bariczbook} \textsc{\'A. Baricz}, Generalized Bessel functions of the first kind, Lecture
 Notes in Mathematics, vol. 1994, Springer, Berlin, 2010.

\bibitem[\bf BP]{bariczpon} \textsc{\'A. Baricz, S. Ponnusamy}, Starlikeness and convexity
of generalized Bessel functions, {\em Integral Transforms Spec.
Funct.} 21(9) (2010) 641--653.

\bibitem[\bf Br1]{brown} \textsc{R.K. Brown}, Univalence of Bessel
functions, {\em Proc. Amer. Math. Soc.} 11(2) (1960) 278--283.

\bibitem[\bf Br2]{brown2} \textsc{R.K. Brown}, Univalent solutions of $W''+pW=0,$ {\em Canad. J. Math.} 14 (1962) 69--78.

\bibitem[\bf Br3]{brown3} \textsc{R.K. Brown}, Univalence of normalized solutions of $W''(z)+p(z)W(z)=0,$
{\em Internat. J. Math. Math. Sci.} 5(3) (1982) 459--483.

\bibitem[\bf HM]{hayden} \textsc{T.L. Hayden, E.P. Merkes}, Chain sequences and univalence, {\em Illinois J. Math.} 8 (1964) 523--528.

\bibitem[\bf IM1]{muldoon} \textsc{M.E.H. Ismail, M.E. Muldoon}, Zeros of combinations of Bessel functions and their derivatives,
{\em Appl. Anal.} 31 (1988) 73--90.

\bibitem[\bf IM2]{ismail} \textsc{M.E.H. Ismail, M.E. Muldoon},
Bounds for the small real and purely imaginary zeros of Bessel and
related functions, {\em Methods Appl. Anal.} 2(1) (1995) 1--21.

\bibitem[\bf Ja]{jack} \textsc{I.S. Jack}, Functions starlike and convex of order $\alpha,$ {\em J. London Math. Soc.} 3(2) (1971) 469--474.

\bibitem[\bf KT]{todd} \textsc{E. Kreyszig, J. Todd}, The radius of univalence of Bessel functions, {\em Illinois J. Math.} 4 (1960) 143--149.

\bibitem[\bf MRS]{merkes} \textsc{E.P. Merkes, M.S. Robertson, W.T. Scott}, On products of starlike functions, {\em Proc. Amer. Math. Soc.} 13 (1962) 960--964.

\bibitem[\bf Ne]{nehari} \textsc{Z. Nehari}, The Schwarzian derivative and schlicht functions, {\em Bull. Amer. Math. Soc.} 55 (1949) 545--551.

\bibitem[\bf Ro]{robertson} \textsc{M.S. Robertson}, Schlicht solutions of $W''+pW=0,$ {\em Trans. Amer. Math. Soc.} 76 (1954) 254--274.

\bibitem[\bf Sz]{szasz1} \textsc{R. Sz\'asz}, On starlikeness of Bessel functions of the
first kind, In: Proceedings of the 8th Joint Conference on Mathematics
and Computer Science, Kom\'arno, Slovakia, 2010, 9pp.

\bibitem[\bf SK]{szasz2} \textsc{R. Sz\'asz, P.A. Kup\'an}, About the univalence of the Bessel functions,
{\em Stud. Univ. Babe\c{s}-Bolyai Math.} 54(1) (2009) 127--132.

\bibitem[\bf Wa]{watson} \textsc{G.N. Watson}, {\em A Treatise on the Theory of Bessel Functions}, 2nd ed., Cambridge University Press, Cambridge, 1944.

\bibitem[\bf Wi]{wilf} \textsc{H.S. Wilf}, The radius of univalence of certain entire functions, {\em Illinois J. Math.} 6 (1962) 242--244.

\end{thebibliography}
\end{document}